\documentclass[11pt]{article} 
\begin{document} 
\newtheorem{prop}{Proposition}[section]
\newtheorem{Def}{Definition}[section] \newtheorem{theorem}{Theorem}[section]
\newtheorem{lemma}{Lemma}[section] \newtheorem{Cor}{Corollary}[section]

\title{\bf Rough solutions of a Schr\"odinger - Benjamin - Ono system}
\author{{\bf Hartmut Pecher}\\
Fachbereich Mathematik und Naturwissenschaften\\
Bergische Universit\"at Wuppertal\\
Gau{\ss}str.  20\\
D-42097 Wuppertal\\
Germany\\
e-mail Hartmut.Pecher@math.uni-wuppertal.de}
\date{}
\maketitle

\begin{abstract}
The Cauchy problem for a coupled Schr\"odinger and Benjamin-Ono system is 
shown to be globally well-posed for a class of data without finite energy. 
The proof uses the I-method introduced by Colliander, Keel, Staffilani, 
Takaoka, and Tao.
\end{abstract}

\renewcommand{\thefootnote}{\fnsymbol{footnote}}
\footnotetext{\hspace{-1.8em}{\it 2000 Mathematics Subject Classification:} 
35Q55, 35Q35 }
\normalsize 
\setcounter{section}{-1}
\section{Introduction}
Consider the following weakly coupled dispersive system 
\begin{eqnarray}
\label{1.1}
i\partial_t u + \partial_x^2 u & = & \alpha uv \\
\label{1.2}
\partial_t v + \nu \partial _x |\partial_x| v & = & \beta \partial_x (|u|^2)
\end{eqnarray}
with Cauchy data
\begin{equation}
u(x,0)  =  u_0(x) \,  , \, v(x,0)  =  v_0(x)
\label{1.3}
\end{equation}
where $ x,t \in {\bf R} $ , $\alpha,\beta, \nu \in {\bf R}$.

This system was introduced by Funakoshi and Oikawa \cite{FO} to model the 
interaction of two fluids described by a short wave term $u: {\bf R} \times 
{\bf 
R} \to {\bf C}$ , which fulfills a Schr\"odinger type equation and a long 
wave 
term $v: {\bf R} \times {\bf R} \to {\bf R}$ , which fulfills a Benjamin - 
Ono 
type equation. Bekiranov, Ogawa and Ponce \cite{BOP} showed local 
well-posedness 
for data $u_0 \in H^s({\bf R}) $ , $ v_0 \in H^{s-\frac{1}{2}}({\bf R}) $ 
and 
$|\nu| \neq 1 $ , $ s \ge 0 \, . $ Because the system satisfies three 
conservation 
laws (cf. (\ref{CL0a}) - (\ref{CL0c}) below) it is not difficult to see that 
this solution exists globally if $\nu > 0$ and $ \frac{\alpha}{\beta} < 0$ 
in the case $ 
s \ge 1$ (finite energy solutions).

In this paper we first show local well-posedness also in the case $|\nu| = 
1$, 
if $ s > 0 $ . Then we use the Fourier restriction norm method and 
especially 
the so-called I-method to show global well-posedness for data with infinite 
energy (and $\nu > 0$ , $ \frac{\alpha}{\beta} < 0$), we assume only $s > 
\frac{1}{3}$. This 
method was introduced by Colliander, Keel, Staffilani, Takaoka, and Tao and 
successfully applied in various situations 
(\cite{CKSTT1},\cite{CKSTT},\cite{CKSTT2},\cite{CKSTT3},\cite{CKSTT4},\cite{
CKSTT5},\cite{CKSTT6},\cite{CST},\cite{P},\cite{P1},\cite{Tz}), in most of 
the cases using a scaling invariance of the problem. Such an invariance is 
also 
very helpful in our case. One introduces for given $0<s<1$ and $N>>1$ the 
Fourier multiplier $\widehat{I_N f}(\xi) := m_N(\xi)\widehat{f}(\xi)$, where 
$m_N$ is a smooth, radially symmetric and nonincreasing function of $|\xi|$ 
and 
$$ m(\xi):=m_N(\xi):= \left\{ \begin{array}{r@{\quad \quad}l} 1 & |\xi| \le 
N \\ 
(\frac{N}{|\xi|})^{1-s} & |\xi| \ge 2N \end{array} \right. $$ 
Then $I=I_N$ is a smoothing operator which maps $H^s({\bf R})$ to $H^1({\bf 
R})$ 
in the sense that
$$ \|u\|_{H^s} \le c \|Iu\|_{H^1} \le c N^{1-s} \|u\|_{H^s} $$
and similarly
$$  \|v\|_{H^{s-{\frac{1}{2}}}} \le c \|Iv\|_{H^{\frac{1}{2}}} \le c N^{1-s} 
\|v\|_{H^{s-\frac{1}{2}}} \, . $$  
One then considers the conserved functionals $L$ and $E$ (cf. (\ref{CL0b}) 
and 
(\ref{CL0c}) below) replacing $u$ and $v$ by $I_N u$ and $I_N v$, so that 
they 
make sense for $u \in H^s$ , $ v \in H^{s-\frac{1}{2}}$, whereas they 
originally 
are only defined for $u \in H^1$ , $ v \in H^{\frac{1}{2}}$ . These modified 
functionals are then shown to be almost conserved in the sense that their 
increment on a local existence interval is bounded by $c N^{-1}$ . One can 
show 
that this is enough to  make the continuation process by reapplying the 
local 
existence theorem uniform, provided $s$ is close enough to 1, namely $s > 
\frac{1}{3}$ .
  
We use the following norms (for $s \in {\bf R} \, , \, -1<b<1$):
$$ \|u\|_{X^{s,b}} := \| \langle \tau + \xi^2 \rangle^b \langle \xi 
\rangle^s 
\widehat{u}(\xi,\tau)\|_{L^2({\bf R}^2)} $$
$$ \|v\|_{Y^{s,b}} := \| \langle \tau + \nu \xi |\xi| \rangle^b \langle \xi 
\rangle^s \widehat{v}(\xi,\tau)\|_{L^2({\bf R}^2)} $$
belonging to the Schr\"odinger and Benjamin - Ono equation, respectively. We 
also need the local in time norm $ \|u\|_{X_{\delta}^{s,b}} := 
\inf_{\psi_{|[0,\delta]}=f} \|\psi\|_{X^{s,b}} $ and similarly 
$\|v\|_{Y^{s,b}_{\delta}}$ .

The standard facts about the Fourier restriction norm method which we use 
without further comments can be found in \cite{BOP}, Chapter 2. The 
Strichartz 
estimates for the homogeneous Schr\"odinger and Benjamin - Ono equation read 
$$ \|e^{it\partial_x^2} u_0\|_{L^6_{xt}} \le c \|u_0\|_{L^2_x} $$ 
and 
$$ \|e^{it\nu \partial_x |\partial_x|}  u_0\|_{L^6_{xt}} \le c 
\|u_0\|_{L^2_x}$$
(cf. \cite{Str} and \cite{KPV1}), which immediately imply $\|u\|_{L^p_{xt}} 
\le c 
\|u\|_{X^{0,b}}$ and $ \|v\|_{L^p_{xt}} \le c \|v\|_{Y^{0,b}}$ for $2 \le p 
\le 
6$ and $b > \frac{1}{2}$. We also use the following bilinear Strichartz type 
estimate 
for the Schr\"odinger equation
\begin{equation}
\label{0}
\|D_x^{1/2}(u_1 \overline{u_2})\|_{L^2_{xt}} \le c \|u_1\|_{X^{0,b}} 
\|u_2\|_{X^{0,b}} \qquad , \qquad b > \frac{1}{2}
\end{equation}
(for a proof cf. e.g. \cite{BOP}, Lemma 3.2).

We denote by $a+$ and $a-$ a number slightly larger and smaller than $a$ , 
respectively.

\section{Local existence}
\begin{prop}
\label{Proposition1}
For $|\nu| = 1$ we have
$$ \|uv\|_{X^{s,a}} \le c \|u\|_{X^{s,b}} \|v\|_{Y^{s-\frac{1}{2},b}} $$
if $ -\frac{1}{2} < a < 0 < \frac{1}{2} < b $ and $ s > 1-2|a|$ 
($\Leftrightarrow |a| > 
\frac{1-s}{2}$) (especially $ s > 0$).
\end{prop}
{\bf Proof:} (along the lines of \cite{BOP}, Lemma 3.4)\\
Assume first $\nu = 1$. We have to prove the following estimate
\begin{eqnarray}
\nonumber
& & \left| \int\int\int\int \frac{\langle\xi\rangle^s g(\sigma,\eta) f(\tau 
- 
\sigma, \xi - \eta) \overline{\phi}(\tau,\xi) \, d\sigma d\eta d\tau 
d\xi}{\langle\eta\rangle^{s-\frac{1}{2}} \langle\xi - \eta\rangle^s \langle 
\tau 
+ \xi^2 \rangle^{|a|} \langle \sigma + \nu \eta |\eta| \rangle^b \langle 
\tau - 
\sigma + (\xi - \eta)^2 \rangle^b} \right| \\ \label{*}
& & \quad \le c \|g\|_{L^2} \|f\|_{L^2} \|\phi\|_{L^2}
\end{eqnarray}
We split $ (\tau,\sigma,\xi,\eta) \in {\bf R}^4 $ into several regions:
\begin{eqnarray*}
A & = & \{ |\eta| < 1 \} \\
B & = & \{ \eta < 0 \, , \, |\xi| \ge \frac{1}{2}|\eta| \, , \, |\eta| \ge 1 
\} 
\\
C & = & \{ \eta < 0 \, , \, |\xi| < \frac{1}{2}|\eta| \, , \, |\eta| \ge 1 
\} \\
D & = & \{ \eta > 0 \, , \, 
|\xi - \eta| \le \frac{1}{2}|\eta| \, , \, |\eta| \ge 1 \} \\
E & = & \{ \eta > 0 \, , \, |\xi - \eta| > \frac{1}{2}|\eta| \, , \, |\eta| 
\ge 
1 \} 
\end{eqnarray*}
Now in $E$ we have
$$ |\nu \eta |\eta| + \eta^2 - 2\xi\eta| = 2|\eta^2-\xi\eta| = 2|\eta||\eta 
- 
\xi| > |\eta|^2 $$ and thus
$$ |\tau + \xi^2| + |\sigma + \nu \eta |\eta|| + |\tau-\sigma+(\xi-\eta)^2| 
\ge 
|\nu \eta |\eta| + \eta^2 -2\xi\eta| > |\eta|^2 \, .$$
According to which of the terms on the l.h.s. is dominant we split $E$ into 
3 
parts:
\begin{eqnarray*}
E_1 & \hspace{-1em} = \hspace{-1em} & E \cap \{ |\tau + \xi^2| \ge |\sigma + 
\nu \eta 
|\eta||,|\tau-\sigma+(\xi-\eta)^2| \,\, , \, |\tau + \xi^2| \ge \frac{1}{3} 
|\eta|^2 \} \\
E_2 & \hspace{-1em} = \hspace{-1em} & E \cap \{|\sigma + \nu \eta |\eta|| 
\ge |\tau + \xi^2|  
,|\tau-\sigma+(\xi-\eta)^2| \,\, , \,  |\sigma + \nu \eta |\eta|| \ge 
\frac{1}{3} 
|\eta|^2 \} \\
E_3 & \hspace{-0.5em} = & \hspace{-0.5em} E \cap\{|\tau-\sigma+(\xi-\eta)^2| 
\ge |\tau + \xi^2|,|\sigma + \nu \eta 
|\eta|| \,\, , \,  |\tau-\sigma+(\xi-\eta)^2| \ge \frac{1}{3} |\eta|^2 \}
\end{eqnarray*}
Define $R_1 = A \cup B \cup D \cup E_1$ , $ R_2 = C \cup E_2$ , $ R_3 = E_3$ 
. 
In order to prove (\ref{*}) in the region $R_1$ it is sufficient to show
\begin{equation}
\label{1}
\left\| \frac{\langle \xi \rangle^s}{\langle \tau + \xi^2 \rangle^{|a|}} 
\left(\int\int \frac{\langle \eta \rangle \chi_{R_1} \, d\sigma 
d\eta}{\langle 
\eta \rangle^{2s} \langle \xi - \eta \rangle^{2s} \langle \sigma + 
\nu\eta|\eta| 
\rangle^{2b} \langle \tau - \sigma + (\xi - \eta)^2\rangle^{2b}} 
\right)^{\frac{1}{2}} \right\|_{L^{\infty}_{\tau}(L_{\xi}^{\infty})} 
\hspace{-0.6em} < \infty
\end{equation} 
Similarly, in order to prove (\ref{*}) in $R_2$ we have to show
\begin{equation}
\label{2}
\left\| \frac{\langle \eta \rangle^{\frac{1}{2}}}{\langle \eta \rangle^s 
\langle 
\sigma + \nu \eta |\eta| \rangle^b} \left( \int\int \frac{\langle \xi 
\rangle^{2s} \chi_{R_2} \, d\tau d\xi}{\langle \xi - \eta \rangle^{2s} 
\langle 
\tau + \xi^2 \rangle^{2|a|} \langle \tau - \sigma + (\xi - \eta)^2 
\rangle^{2b}} 
\right)^{\frac{1}{2}} \right\|_{L_{\sigma}^{\infty}(L_{\eta}^{\infty})} 
\hspace{-0.4em} < \infty
\end{equation} 
Finally, in order to prove (\ref{*}) in $R_3$ we use the transformed region
$$ \widetilde{R_3} = \{ (\rho,\sigma,\zeta,\eta) \in {\bf R}^4 : |\eta|\ge 1 
\, , \, 
|\rho - \zeta^2| \ge \frac{1}{3} |\eta|^2 \} \, , $$
where $\rho := \sigma - \tau$ , $ \zeta := \eta - \xi $. We have to show
\begin{equation}
\label{3}
\left\|\frac{1}{\langle \zeta \rangle^s \langle \rho - \zeta^2 \rangle^b} 
\left(\int\int \frac{\langle \eta - \zeta \rangle^{2s} 
\chi_{\widetilde{R_3}} \, 
d\sigma d\eta}{\langle \eta \rangle^{2s-1} \langle \sigma + \nu \eta |\eta| 
\rangle^{2b} \langle \sigma - \rho + (\eta - \zeta)^2\rangle^{2|a|}} 
\right)^{\frac{1}{2}} \right\|_{L^{\infty}_{\rho}(L^{\infty}_{\zeta})} < 
\infty \, 
.
\end{equation} 
We start to prove (\ref{1}). In the regions $A,B$ and $D$ we use the 
estimate 
$\langle \xi \rangle \le \langle \eta \rangle \langle \xi - \eta \rangle$ so 
that it suffices to show
$$ \left\| \int\int \frac{\langle \eta \rangle \, d\sigma d\eta}{\langle 
\sigma 
+ \nu \eta |\eta|\rangle^{2b} \langle \tau - \sigma +(\xi - 
\eta)^2\rangle^{2b}} 
\right\|_{L^{\infty}_{\tau}(L^{\infty}_{\xi})} < \infty \, . $$
Performing the $\sigma$-integration we get by \cite{BOP}, Lemma 2.5 (2.11):
$$ \int\int \frac{\langle \eta \rangle \, d\sigma d\eta}{\langle \sigma + 
\nu 
\eta |\eta|\rangle^{2b} \langle \tau - \sigma +(\xi - \eta)^2\rangle^{2b}} 
\le c 
\int \frac{\langle \eta \rangle \, d\eta}{\langle \tau + \xi^2 + \nu \eta 
|\eta| 
+ \eta^2 - 2\xi\eta \rangle^{2b}} \, .$$
This is trivially bounded in the region $A$, whereas in region $B$ we 
substitute 
$\tau + \xi^2 + \nu \eta |\eta| + \eta^2 -2\xi\eta = \tau + \xi^2 -2\xi\eta 
=: 
\eta' $, so that $ \frac{d\eta'}{d\eta} = -2\xi$, and we get the bound using 
$|\xi| \ge \frac{1}{2}|\eta|$ and $|\eta| \ge 1$:
$$ c \int \frac{\langle \eta \rangle \, d\eta'}{\langle\eta'\rangle^{2b} 
|\xi|} 
\le c \int \frac{d\eta'}{\langle\eta'\rangle^{2b}} < \infty $$
for $b > \frac{1}{2}$.\\
In the region $D$ we have $ \tau + \xi^2 + \nu \eta |\eta| + \eta^2 - 
2\xi\eta = 
\tau + \xi^2 + 2\eta^2 - 2\xi\eta =: \eta' $, so that $ 
|\frac{d\eta'}{d\eta}| = 
|4\eta-2\xi| =|2(2\eta-\xi)|=2|\eta+(\eta-\xi)| \ge 2(|\eta|-|\eta-\xi|) \ge 
|\eta| $, because $|\xi-\eta| \le \frac{1}{2}|\eta|$. Thus we get the bound
$$ c \int \frac{\langle \eta\rangle \, d\eta'}{\langle\eta'\rangle^{2b} 
|\eta|} 
< \infty $$
for $ b > \frac{1}{2} $ and $ |\eta| \ge 1 $ . \\
It remains to prove (\ref{1}) in the region $E_1$. First we consider the 
case 
$0<s<\frac{1}{2}$. This implies $|a|> \frac{1-s}{2} > \frac{1}{4}$. We use 
again 
the estimate $\langle \xi \rangle \le \langle \eta \rangle \langle \xi - 
\eta 
\rangle $ , so that it suffices to show
$$ \left\| \frac{1}{\langle \tau + \xi^2 \rangle^{|a|}} \left(\int\int 
\frac{\langle\eta\rangle \, d\sigma d\eta}{\langle \sigma + \nu \eta 
|\eta|\rangle^{2b} \langle\tau - \sigma + 
(\xi-\eta)^2\rangle^{2b}}\right)^{\frac{1}{2}}\right\|_{L^{\infty}_{\tau} 
(L^{\infty}_{\xi})} < \infty \, . $$ 
Performing the $\sigma$-integration as above it remains to bound
\begin{eqnarray*}
\lefteqn{\frac{1}{\langle \tau + \xi^2 \rangle^{|a|}} \left(\int 
\frac{\langle 
\eta \rangle \, d\eta}{\langle \tau + \xi^2 + \nu \eta |\eta| + \eta^2 - 
2\xi\eta \rangle^{2b}} \right)^{\frac{1}{2}}} \\
& \le & c \left( \int \frac{d\eta}{\langle \tau + \xi^2 + 2\eta^2-2\xi\eta 
\rangle^{2b}}\right)^{\frac{1}{2}} \le c 
\end{eqnarray*}
using $\langle \tau + \xi^2 \rangle^{|a|} \ge c \langle \eta \rangle^{2|a|} 
\ge 
c \langle \eta \rangle^{\frac{1}{2}}$ . \\
Next we consider the case $s \ge \frac{1}{2} $ in the region $E_1$ . First 
of 
all, consider the subregion $|\xi| \ge \frac{3}{2} |\eta|$ . In this case we 
get 
the following bound for (\ref{1}) performing the $\sigma$-integration:
\begin{eqnarray*}
\lefteqn{c \left(\int\int \frac{ \langle \eta \rangle \, d\sigma 
d\eta}{\langle 
\eta \rangle^{2s} \langle \sigma + \nu \eta |\eta| \rangle^{2b} \langle \tau 
- 
\sigma +(\xi - \eta)^2 \rangle^{2b}} \right)^{\frac{1}{2}}} \\ & \le & c 
\left( 
\int \frac{\langle \eta \rangle \, d\eta} {\langle \eta \rangle^{2s} \langle 
\tau + \xi^2 + \nu \eta |\eta| + \eta^2 - 2\xi\eta \rangle^{2b}} 
\right)^{\frac{1}{2}} \\
& \le & c \left( \int \frac{d\eta}{\langle \tau + \xi^2 + 2\eta^2 - 2\xi\eta 
\rangle^{2b}} \right)^{\frac{1}{2}} \le c \, .
\end{eqnarray*}
In the subregion $|\xi| \le \frac{3}{2}|\eta|$ we get the following bound 
for 
(\ref{1}) performing the $\sigma$-integration:
\begin{eqnarray*}
\lefteqn{c \left( \int\int \frac{\langle \eta \rangle \, d\sigma 
d\eta}{\langle 
\sigma + \nu \eta |\eta|\rangle^{2b} \langle \tau - \sigma +(\xi - 
\eta)^2\rangle^{2b}} \right)^{\frac{1}{2}}} \\
& \le & c \left( \int \frac{ \langle \eta \rangle \, d\eta}{\langle \tau + 
\xi^2 
+ \nu \eta |\eta| + \eta^2 - 2\xi\eta \rangle^{2b}} \right)^{\frac{1}{2}} \, 
. 
\end{eqnarray*}
Substituting $ \tau +\xi^2 + \nu \eta |\eta| + \eta^2 - 2\xi\eta = \tau + 
\xi^2 
+ 2\eta^2 - 2\xi\eta =: \eta' $ we have $|\frac{d\eta'}{d\eta}| = 
|4\eta-2\xi| 
= 2(|2\eta|-|\xi|) \ge 2(|2\eta| - \frac{3}{2}|\eta|) = |\eta| \sim \langle 
\eta \rangle $, and we get the bound $ c \left(\int  \langle \eta' 
\rangle^{-2b} 
d\eta' \right)^{\frac{1}{2}} < \infty $ . \\
Next we have to prove (\ref{2}). In the region $C$ it suffices to show
\begin{equation}
\label{**}
\left\| \frac{\langle \eta \rangle^{\frac{1}{2}}}{\langle \eta \rangle^{2s}} 
\left( \int\int \frac{\langle \xi \rangle^{2s}\chi_C \, d\tau d\xi}{\langle 
\tau 
+ \xi^2 \rangle^{2|a|} \langle \tau - \sigma +(\xi-\eta)^2 \rangle^{2b}} 
\right)^{\frac{1}{2}} \right\|_{L^{\infty}_{\sigma}L^{\infty}_{\eta}} < 
\infty \, .
\end{equation} 
The integration with respect to $\tau$ gives by \cite{BOP}, Lemma 2.5 (2.11) 
and 
H\"older:
\begin{eqnarray*}
\int\int \frac{\langle\xi\rangle^{2s} \chi_C \, d\tau d\xi}{\langle \tau + 
\xi^2 
\rangle^{2|a|} \langle \tau - \sigma + (\xi - \eta)^2 \rangle^{2b}}  \le c 
\int 
\frac{\langle \xi \rangle^{2s} \chi_C \, d\xi}{\langle \sigma - \eta^2 + 
2\xi\eta \rangle^{2|a|}} \\
 \le c \left(\int_{|\xi| \le \frac{1}{2} |\eta|} \langle \xi 
\rangle^{2s\widehat{p}} \, d\xi \right)^{\frac{1}{\widehat{p}}} \left( \int 
\frac{d\xi}{\langle \sigma - \eta^2 + 2\xi\eta \rangle^{2|a|\widehat{q}}} 
d\xi 
\right)^{\frac{1}{\widehat{q}}} \, . 
\end{eqnarray*}
Choosing $ \frac{1}{\widehat{q}} = 2|a|- $ , $ \frac{1}{\widehat{p}} = 
1+2|a|+ $ 
and substituting $\eta' = \sigma - \eta^2 + 2\xi\eta $ we get the bound
\begin{eqnarray*}
c \langle \eta \rangle^{\frac{2s\widehat{p}+1}{\widehat{p}}} \left(\int 
\frac{d\eta'}{\langle \eta' \rangle^{1+} |\eta|} 
\right)^{\frac{1}{\widehat{q}}} 
& \le & c \langle \eta 
\rangle^{2s+\frac{1}{\widehat{p}}-\frac{1}{\widehat{q}}} 
\\
& = & c \langle \eta \rangle^{2s+1-\frac{2}{\widehat{q}}} = c \langle \eta 
\rangle^{2s+1-4|a|+} \, .
\end{eqnarray*}
Thus we get the following bound for (\ref{**}):
$$ \langle \eta \rangle^{\frac{1}{2}-2s+s+\frac{1}{2}-2|a|+} = \langle \eta 
\rangle^{1-s-2|a|+} \le c $$
because $ |a| > \frac{1-s}{2}$ . \\
Next we prove (\ref{2}) in the region $E_2$ . Using the estimate $ \langle 
\xi 
\rangle \le \langle \eta \rangle \langle \xi-\eta \rangle $ and performing 
the 
$\tau$-integration it suffices to get a bound on $E_2$ for
$$ \frac{\langle \eta \rangle^{\frac{1}{2}}}{\langle \sigma + \nu \eta 
|\eta| 
\rangle^b} \left( \int \frac{d\xi}{\langle \sigma - \eta^2 + 2\eta\xi 
\rangle^{2|a|}} \right)^{\frac{1}{2}} \, .$$
Substitution of $ \eta' := \sigma - \eta^2 +2\eta\xi $ and using the 
definition 
of $E_2$ we get
\begin{eqnarray*}
|\eta'|=|\sigma - \eta^2 + 2\eta\xi| &  = & |(\tau + 
\xi^2)-(\tau-\sigma+(\xi-\eta)^2)| \\ & \le & |\tau + \xi^2| 
+|\tau-\sigma+(\xi-\eta)^2| \le 2|\sigma+\nu\eta|\eta|| 
\end{eqnarray*}
and thus
\begin{eqnarray*}
& & \frac{\langle \eta \rangle^{\frac{1}{2}}}{\langle \sigma + \nu \eta 
|\eta|\rangle^b} \left(\int \frac{d\xi}{\langle \sigma - \eta^2 + 2\eta\xi 
\rangle^{2|a|}} \right)^{\frac{1}{2}} \\
 & & \le  c \frac{\langle \eta \rangle^{\frac{1}{2}}}{\langle \sigma + \nu 
\eta 
|\eta|\rangle^b}
 \left(\int_{|\eta'|  \le  2|\sigma+\nu\eta|\eta||} \frac{d\eta'}{\langle 
\eta'\rangle^{2|a|}|\eta|} \right)^{\frac{1}{2}} 
 \le  c \frac{\langle\sigma + 
\nu\eta|\eta|\rangle^{\frac{1}{2}-|a|}}{\langle 
\sigma + \nu\eta|\eta|\rangle^b} \le c \, .
\end{eqnarray*}
It remains to prove (\ref{3}) in the region $\widetilde{R_3}$. Using the 
estimate 
$\langle \eta - \zeta\rangle \le \langle \zeta \rangle \langle \eta \rangle$ 
and 
performing the $\sigma$-integration it is enough to give the following bound 
in 
$\widetilde{R_3}$ :
$$ \frac{1}{\langle \rho - \zeta^2 \rangle^b} \left( \int \frac{\langle 
\eta\rangle \, d\eta}{\langle \rho - \zeta^2 + \nu \eta |\eta| - \eta^2 + 
2\zeta 
\eta \rangle^{2|a|}} \right)^{\frac{1}{2}} \le c \left(\int 
\frac{d\eta}{\langle 
\eta \rangle^{4b-1}} \right)^{\frac{1}{2}} \le c \, .$$
The case $\nu = -1$ can be treated in the same way by replacing $\eta<0$ by 
$\eta>0$ in the regions $B$ and $C$ and $\eta>0$ by $\eta<0$ in $D$ and $E$.
\bigskip

Moreover the following estimates for the nonlinearities are true (cf. 
\cite{BOP}, Cor. 3.3 and Lemma 3.4).
\begin{prop}
\label{Proposition2}
\begin{enumerate}
\item For arbitrary $\nu$ and $s\ge 0$ , $b > \frac{1}{2}$ :
$$ \| (|u|^2)_x \|_{Y^{s-\frac{1}{2},0}} \le c \|u\|^2_{X^{s,b}} \, . $$
\item For $|\nu| \neq 1$ and $s\ge 0$ , $ b > \frac{1}{2}$ :
$$ \|uv\|_{X^{s,-\frac{1}{4}}} \le c \|u\|_{X^{s,b}} 
\|v\|_{Y^{s-\frac{1}{2},b}} \, .
$$
\end{enumerate}
\end{prop}

These propositions imply by standard arguments the following local existence 
result.
\begin{theorem}
\label{Theorem1}
Let $s>0$ in the case $|\nu|=1$, and $s \ge 0$ in the case $|\nu| \neq 1$ . 
For 
any $(u_0,v_0) \in H^s({\bf R}) \times H^{s-\frac{1}{2}}({\bf R})$ there 
exists $ 
b > \frac{1}{2} $ and $\, \delta = 
\delta(\|u_0\|_{H^s},\|v_0\|_{H^{s-\frac{1}{2}}})$ $ > 0 $ such that the 
Cauchy 
problem (\ref{1.1}),(\ref{1.2}),(\ref{1.3}) has a unique solution $(u,v) \in 
X^{s,b}_{\delta} \times Y^{s-\frac{1}{2},b}_{\delta}$ and $(u,v) \in 
C^0([0,\delta],H^s \times H^{s-\frac{1}{2}})$ .
\end{theorem}

Applying the operator $I$ to the system (\ref{1.1}),(\ref{1.2}),(\ref{1.3}) 
we 
get the problem
\begin{eqnarray}
\label{2.1}
i I\partial_t u + I\partial^2_xu & = & \alpha I(uv) \\
\label{2.2}
I\partial_t v + \nu I(\partial_x |\partial_x| u) & = & \beta I \partial_x 
(|u|^2) \\
\label{2.3}
Iu(0) = Iu_0 & , & Iv(0) = Iv_0 \, .
\end{eqnarray}
For this system the following modified local existence result holds:
\begin{prop}
\label{Proposition3}
Assume $1 \ge s>0$, if $|\nu| =1$, and $s \ge 0$ otherwise. For any 
$(u_0,v_0) \in 
H^s 
\times H^{s-\frac{1}{2}}$ there exists $\delta \le 1$ and $ \delta \sim 
(\|Iu_0\|_{H^1} + \|Iv_0\|_{H^{\frac{1}{2}}})^{-\frac{2}{s}-} \, ,$ if 
$|\nu| = 
1$ , and $\delta \sim (\|Iu_0\|_{H^1} + \|Iv_0\|_{H^{\frac{1}{2}}})^{-4-}$ , 
if 
$|\nu| \neq 1$ , such that the system (\ref{2.1}),(\ref{2.2}),(\ref{2.3}) 
has a 
unique local solution in the time interval $[0,\delta]$ with the property
$$ \|Iu\|_{X^{1,b}_{\delta}} + \|Iv\|_{Y^{\frac{1}{2},b}_{\delta}} \le 
\widehat{c} (\|Iu_0\|_{H^1} + \|Iv_0\|_{H^{\frac{1}{2}}}) \, ,  $$
where $ b = \frac{1}{2}+$ .
\end{prop}
{\bf Proof:} We construct a fixed point of the mapping 
$S=(\tilde{S_0},\tilde{S_1})$ induced by the integral equations belonging to 
the 
system (\ref{2.1}),(\ref{2.2}),(\ref{2.3}):
\begin{eqnarray*}
\tilde{S_0}(Iu,Iv)(t) & := & e^{it\partial_x^2} Iu_0 - i \int_0^t 
e^{i(t-s)\partial_x^2} \alpha I(u(s)v(s)) \, ds \\
\tilde{S_1}(Iu,Iv)(t) & := & e^{-\nu t\partial_x |\partial_x|} Iv_0 + 
\int_0^t 
e^{-\nu(t-s)\partial_x |\partial_x|} \beta I(|u(s)|^2)_x \, ds \, .
\end{eqnarray*}
The estimates for the nonlinearities in (\ref{Proposition1}) and 
(\ref{Proposition2}) carry over to corresponding estimates including the 
$I$-operators by the interpolation lemma of \cite{CKSTT5}, namely
\begin{eqnarray*}
\|I(|u|^2)_x\|_{Y^{\frac{1}{2},0}} & \le & c \|Iu\|_{X^{1,b}}^2 \\
\|I(uv)\|_{X^{1,-|a|}} & \le & c \|Iu\|_{X^{1,b}} \|Iv\|_{Y^{\frac{1}{2},b}}
\end{eqnarray*}
(with $|a|=\frac{1-s}{2}+$ , if $|\nu|=1$, and $|a|=\frac{1}{4}$ otherwise). 
This implies
\begin{eqnarray*}
\|\tilde{S_0}(Iu,Iv)\|_{X^{1,b}_{\delta}} & \le & c \|Iu_0\|_{H^1} + c 
|\alpha| 
\|Iu\|_{X^{1,b}_{\delta}} \|Iv\|_{Y^{\frac{1}{2},b}_{\delta}} 
\delta^{\frac{1}{2}-|a|-} \\
\|\tilde{S_1}(Iu,Iv)\|_{Y^{\frac{1}{2},b}_{\delta}} & \le & c 
\|Iv_0\|_{H^{\frac{1}{2}}} + c |\beta| \|Iu\|_{X^{1,b}_{\delta}}^2 
\delta^{\frac{1}{2}-} \, .
\end{eqnarray*}   
This gives the desired bounds, provided
$$ c \delta^{\frac{1}{2}-|a|-} (\|Iu_0\|_{H^1} + \|Iv_0\|_{H^{\frac{1}{2}}}) 
< 1 
\, . $$

\section{Conservation laws}
Our system has the following conserved quantities:
\begin{eqnarray}
\label{CL0a}
M & := & \|u\| \\
\label{CL0b}
L(u,v) & := & -\frac{\alpha}{2\beta} \|v\|^2 - \Im \int u \overline{u_x} \, 
dx 
\\
\label{CL0c}
E(u,v) & := & \|u_x\|^2 - \frac{\alpha \nu}{2\beta} \|D_x^{\frac{1}{2}} 
v\|^2 
+\alpha \int v |u|^2 \, dx \, .
\end{eqnarray}
From now on, we assume $\nu > 0$ and $\frac{\alpha}{\beta} < 0$ .\\
Then $L$ and $E$ are controlled by $\|u\|_{H^1}$ and 
$\|v\|_{H^{\frac{1}{2}}}$, 
and vice versa, as one concludes as follows:
\begin{equation}
\label{CL1}
|L(u,v)| \le c \|v\|^2 + M \|u_x\|
\end{equation}
and
\begin{equation}
\label{CL2}
\|v\|^2 \le c(|L| + M \|u_x\|) \, .
\end{equation}
Concerning $E$ we have by Gagliardo - Nirenberg
\begin{eqnarray*}
\int |vu^2| \, dx & \le & \|v\| \|u\|^{\frac{3}{2}} \|u_x\|^{\frac{1}{2}} 
\le 
c(|L|^{\frac{1}{2}} M^{\frac{3}{2}} \|u_x\|^{\frac{1}{2}} + M^2 \|u_x\|) \\
& \le & c(|L|^{\frac{2}{3}} M^2 + M^4) + \epsilon \|u_x\|^2 \\
& \le & c(|L|^{\frac{4}{3}} + M^4) + \epsilon \|u_x\|^2
\end{eqnarray*}
and thus
$$ \|u_x\|^2 + |\frac{\alpha \nu}{\beta}| \|D_x^{\frac{1}{2}} v\|^2 \le  |E| 
+ 
c(|L|^{\frac{4}{3}} + M^4) + \epsilon \|u_x\|^2 \, ,$$
consequently
\begin{equation}
\label{CL3}
\|u_x\|^2 + \|D_x^{\frac{1}{2}} v\|^2 \le c(|E| + |L|^{\frac{4}{3}} + M^4) 
\, 
.
\end{equation}
Similarly
\begin{equation}
\label{CL4}
|E| \le c(\|u_x\|^2 + \|D_x^{\frac{1}{2}} v\|^2 + |L|^{\frac{4}{3}} + M^4) 
\, 
.
\end{equation}
From (\ref{CL1}) and (\ref{CL4}) we get
\begin{eqnarray}
\nonumber
|E| & \le &  c (\|u_x\|^2 + \|D_x^{\frac{1}{2}} v\|^2 + \|v\|^{\frac{8}{3}} 
+ 
M^{\frac{4}{3}} \|u_x\|^{\frac{4}{3}} + M^4) \\
\label{CL5}
& \le & c (\|u_x\|^2 + \|D_x^{\frac{1}{2}} v\|^2 + \|v\|^{\frac{8}{3}} + 
M^4) 
\, .
\end{eqnarray}
From (\ref{CL2}) and (\ref{CL3}) we have
\begin{equation}
\|v\|^2  \le  c (|L| + M(|E|^{\frac{1}{2}} + |L|^{\frac{2}{3}} + M^2)) \le  
c (|L| + M |E|^{\frac{1}{2}} + M^3 +1)\, .
 \label{CL6}
\end{equation}
Finally, from (\ref{CL3}) and (\ref{CL6}) we arrive at
\begin{equation}
\label{CL7}
\|u\|_{H^1}^2 + \|v\|_{H^{\frac{1}{2}}}^2 \le c (|E| + |L|^{\frac{4}{3}} + 
M^4 
+1) \, .
\end{equation} 

These estimates imply a-priori-bounds for the $H^1$-norm of $u$ and the 
$H^{\frac{1}{2}}$-norm of $v$ for data with finite energy $E$, finite $L$ 
and 
finite $\|u_0\|$ . This is the case for $H^1$-data $u_0$ and 
$H^{\frac{1}{2}}$-data $v_0$. Thus our local existence result implies
\begin{theorem}
\label{Theorem2}
For data $(u_0,v_0) \in H^1 \times H^{\frac{1}{2}}$ and $\nu >0$ , $ 
\frac{\alpha}{\beta} < 0$ there exists $ b > \frac{1}{2} $ such that 
(\ref{1.1}),(\ref{1.2}),(\ref{1.3}) has a unique global solution $(u,v) \in 
X^{1,b} \times Y^{\frac{1}{2},b}$ with $(u,v) \in C^0({\bf R^+},H^1 \times 
H^{\frac{1}{2}})$ .
\end{theorem}

In order to get a corresponding result for less regular data we consider the 
modified functionals $E(Iu,Iv)$ and $L(Iu,Iv)$ .

Using the modified system (\ref{2.1}), (\ref{2.2}) an elementary calculation 
shows
\begin{eqnarray}
\nonumber
 \frac{d}{dt} E(Iu,Iv) 
& = & \alpha \nu \langle I(|u|^2)-|Iu|^2,D_x Iv_x \rangle + \alpha \beta 
\langle 
I(|u|^2)_x - (|Iu|^2)_x,|Iu|^2  \rangle \\ 
 \nonumber
& & -2\alpha^2 \Im \langle IvIu,I(vu)-IvIu \rangle - 2\alpha \Im \langle 
Iu_x,I(vu)_x - (IvIu)_x \rangle \\
& =: & \sum_{j=1}^4 I_j
\label{CL8}
\end{eqnarray}
and
\begin{equation}
\label{CL9}
\frac{d}{dt} L(Iu,Iv) = -\alpha(\langle Iv,(I(|u|^2)-(|Iu|^2)_x \rangle + 2 
\Re 
\langle I(vu)-IvIu,Iu_x \rangle ) \, .
\end{equation}

\section{Almost conservation}
\begin{prop}
\label{Proposition4}
If $(u,v)$ is a solution of (\ref{1.1}),(\ref{1.2}),(\ref{1.3}) in 
$[0,\delta]$ 
in the sense of Proposition \ref{Proposition3}, then the following estimate 
holds for $N \ge 1$ and $ s \ge \frac{1}{4} $ :
\begin{eqnarray*}
\lefteqn{ |E(Iu(\delta),Iv(\delta)) - E(Iu(0),Iv(0))| + 
|L(Iu(\delta),Iv(\delta)) - L(Iu(0),Iv(0))| } \\
& \le & c N^{-1} \left(\|Iu\|^2_{X^{1,\frac{1}{2}+}_{\delta}} 
\|Iv\|_{Y^{\frac{1}{2},\frac{1}{2}+}_{\delta}} + 
\|Iu\|^4_{X^{1,\frac{1}{2}+}_{\delta}} + 
\|Iu\|^2_{X^{1,\frac{1}{2}+}_{\delta}} 
\|Iv\|_{Y^{\frac{1}{2},\frac{1}{2}+}_{\delta}}^2\right) \, .
\end{eqnarray*}
\end{prop}
{\bf Proof:} Integrating (\ref{CL8}) over $t \in [0,\delta]$ we have to 
estimate 
the various terms on the right hand side. We assume w.l.o.g. the Fourier 
transforms of all the functions to be nonnegative. We drop $\delta$ from the 
notation $X^{s,b}_{\delta}$ and $Y^{s,b}_{\delta}$ .\\
\underline{Estimate of $I_1$:} We have to show
\begin{eqnarray*}
& & \int_0^{\delta} \int_* \left|\frac{m(\xi_1 + \xi_2) - 
m(\xi_1)m(\xi_2)}{m(\xi_1)m(\xi_2)}\right| |\xi_1 + \xi_2| 
\widehat{u_1}(\xi_1,t) \widehat{\overline{u_2}}(\xi_2,t) |\xi_3| 
\widehat{v}(\xi_3,t) \, d\xi dt \\
& & \le  c N^{-1} \|u_1\|_{X^{1,\frac{1}{2}+}} \|u_2\|_{X^{1,\frac{1}{2}+}} 
\|v\|_{Y^{\frac{1}{2},\frac{1}{2}+}} \, .
\end{eqnarray*}
Here and in the sequel * denotes integration over the set $ \sum \xi_i = 0 $ 
. 
We may assume $|\xi_1| \ge N$ or $|\xi_2| \ge N$, because otherwise the 
multiplier term vanishes, and also the two largest frequencies are 
equivalent.\\
\underline{Case 1:} $|\xi_1| << |\xi_2| \sim |\xi_3| $ , $ |\xi_1| \le N $ , 
$ 
|\xi_2| \ge N $ . \\
Using the mean value theorem the multiplier term is estimated  by
$$ \left| \frac{m(\xi_1+\xi_2)-m(\xi_2)}{m(\xi_2)}\right| \le c \left| 
\frac{(\nabla m)(\xi_2)\xi_1}{m(\xi_2)} \right| \le c 
\frac{|\xi_1|}{|\xi_2|} 
\le c \frac{|\xi_1|}{N} \, .$$
Thus by use of (\ref{0}) the integral is bounded by
\begin{eqnarray*}
& & \frac{c}{N} \int_0^{\delta}\int_* |\xi_1 + \xi_2|^{\frac{1}{2}} |\xi_1| 
\widehat{u_1}(\xi_1,t) |\xi_2| \widehat{\overline{u_2}}(\xi_2,t) 
|\xi_3|^{\frac{1}{2}} \widehat{v}(\xi_3,t) \, d\xi dt \\
& & \le \frac{c}{N} \|D_x^{\frac{1}{2}}(D_x u_1 D_x 
\overline{u_2})\|_{L^2_{xt}} 
\|D_x^{\frac{1}{2}} v\|_{L^2_{xt}} \\
& & \le \frac{c}{N} \|u_1\|_{X^{1,\frac{1}{2}+}} 
\|u_2\|_{X^{1,\frac{1}{2}+}} 
\|v\|_{Y^{\frac{1}{2},\frac{1}{2}+}} \, .
\end{eqnarray*}
\underline{Case 2:} $|\xi_1| << |\xi_2| \sim |\xi_3|$ , $ |\xi_1|,|\xi_2| 
\ge N$ 
. \\
The multiplier is bounded by $ \frac{c}{m(\xi_1)} \le c \frac{|\xi_1|}{N} $ 
. 
Thus we can conclude as in Case 1. \\
\underline{Case 3:} $|\xi_1| \sim |\xi_2| \ge cN $ , $ |\xi_1 + \xi_2| \le 
2N $ 
($\Longrightarrow |\xi_3| \le c|\xi_1|,c|\xi_2|$). \\
The multiplier is bounded by $ \frac{c}{m(\xi_1)m(\xi_2)} \le c 
\frac{|\xi_1| 
|\xi_2|}{N^2} $ . Thus we get the following bound for the integral using 
(\ref{0}):
\begin{eqnarray*}
& & \frac{c}{N^2} \int_0^{\delta}\int_* N |\xi_1| \widehat{u_1}(\xi_1,t) 
|\xi_2| 
\widehat{\overline{u_2}}(\xi_2,t) |\xi_1 + \xi_2|^{\frac{1}{2}} 
|\xi_3|^{\frac{1}{2}} \widehat{v}(\xi_3,t) \, d\xi dt \\
& & \le \frac{c}{N} \|D_x^{\frac{1}{2}} (D_x u_1 D_x 
\overline{u_2})\|_{L^2_{xt}} \|D_x^{\frac{1}{2}} v\|_{L^2_{xt}} \\
& & \le \frac{c}{N}  \|u_1\|_{X^{1,\frac{1}{2}+}} 
\|u_2\|_{X^{1,\frac{1}{2}+}} 
\|v\|_{Y^{\frac{1}{2},\frac{1}{2}+}} \, .
\end{eqnarray*}
\underline{Case 4:} $|\xi_1| \sim |\xi_2| \ge cN $ , $ |\xi_1 + \xi_2| \ge 
2N $ . 
\\
The multiplier is bounded by
$$ \frac{m(\xi_1 + \xi_2)}{m(\xi_1)m(\xi_2)} + 1 \le c \frac{|\xi_1|^{1-s} 
|\xi_2|^{1-s} N^{1-s}}{N^{1-s} N^{1-s} |\xi_1 + \xi_2|^{1-s}} = c 
\frac{|\xi_1|^{1-s} |\xi_2|^{1-s}}{N^{1-s} |\xi_1 + \xi_2|^{1-s}} \, , $$
which gives the integral bound
\begin{eqnarray*}
& & \frac{c}{N^{1-s}} \int_0^{\delta} \int_* |\xi_1 + \xi_2|^s |\xi_1|^{1-s} 
\widehat{u_1}(\xi_1,t) |\xi_2|^{1-s} \widehat{\overline{u_2}}(\xi_2,t) 
|\xi_3| 
\widehat{v}(\xi_3,t) \, d\xi dt \\
& & \le \frac{c}{N^{1-s}} \int_0^{\delta} \int_* |\xi_1| 
\widehat{u_1}(\xi_1,t) 
\frac{|\xi_2|}{N^s} \widehat{\overline{u_2}}(\xi_2,t) |\xi_1 + 
\xi_2|^{\frac{1}{2}} |\xi_3|^{\frac{1}{2}} \widehat{v}(\xi_3,t) \, d\xi dt 
\\
& & \le \frac{c}{N} \|D_x^{\frac{1}{2}} (D_x u_1 D_x 
\overline{u_2})\|_{L^2_{xt}} \|D_x^{\frac{1}{2}} v\|_{L^2_{xt}} \\
& & \le \frac{c}{N} \|u_1\|_{X^{1,\frac{1}{2}+}} 
\|u_2\|_{X^{1,\frac{1}{2}+}} 
\|v\|_{Y^{\frac{1}{2},\frac{1}{2}+}} \, .
\end{eqnarray*}
\underline{Estimate of $I_4$:} It is sufficient to show
\begin{eqnarray*}
& & \int_0^{\delta} \int_* \left| \frac{m(\xi_1+\xi_2) - 
m(\xi_1)m(\xi_2)}{m(\xi_1)m(\xi_2)} \right| |\xi_1 + \xi_2| 
\widehat{v}(\xi_1,t) 
\widehat{u_2}(\xi_2,t) |\xi_3| \widehat{\overline{u_3}}(\xi_3,t) \, d\xi dt 
\\
& & \le c N^{-1} \|u_2\|_{X^{1,\frac{1}{2}+}} \|u_3\|_{X^{1,\frac{1}{2}+}} 
\|v\|_{Y^{\frac{1}{2},\frac{1}{2}+}} \, .
\end{eqnarray*}
\underline{Case 1:} $ |\xi_1| << |\xi_2| \sim |\xi_3| \, , \, |\xi_2| \ge N 
$ 
($\Longrightarrow |\xi_1 + \xi_2| \sim |\xi_2|$). \\
If $|\xi_1| \le N$, the multiplier is bounded by
$$ \left| \frac{m(\xi_1 + \xi_2)-m(\xi_2)}{m(\xi_2)} \right| \le c \left| 
\frac{(\nabla m)(\xi_2) \xi_1}{m(\xi_2)} \right| \le c 
\frac{|\xi_1|}{|\xi_2|} 
\le c \frac{|\xi_1|}{N} $$
and, if $|\xi_1| \ge N$, we have the bound $ \frac{c}{m(\xi_1)} \le c 
\frac{|\xi_1|}{N} $ , so that the integral is bounded by
\begin{eqnarray*}
\lefteqn{ \frac{c}{N} \int_0^{\delta} \int_* |\xi_1| |\xi_1 + \xi_2| 
\widehat{v}(\xi_1,t) \widehat{u_2}(\xi_2,t) |\xi_3| 
\widehat{\overline{u_3}}(\xi_3,t) \, d\xi dt } \\
& \le & \frac{c}{N} \int_0^{\delta} \int_* |\xi_1|^{\frac{1}{2}} |\xi_2 + 
\xi_3|^{\frac{1}{2}} \widehat{v}(\xi_1,t) |\xi_2| \widehat{u_2}(\xi_2,t) 
|\xi_3| 
\widehat{\overline{u_3}}(\xi_3,t) \, d\xi dt \\
& \le & \frac{c}{N} \|D_x^{\frac{1}{2}} v\|_{L^2_{xt}} 
\|D_x^{\frac{1}{2}}(D_x 
u_2 D_x \overline{u_3})\|_{L^2_{xt}} \\
& \le & \frac{c}{N} \|u_2\|_{X^{1,\frac{1}{2}+}} 
\|u_3\|_{X^{1,\frac{1}{2}+}} 
\|v\|_{Y^{\frac{1}{2},\frac{1}{2}+}} \, .
\end{eqnarray*}
\underline{Case 2:} $ |\xi_1| >> |\xi_2| $ , $ |\xi_3| \sim |\xi_1| $ , $ 
|\xi_1| \ge N $ ($\Longrightarrow |\xi_1 + \xi_2| \sim |\xi_1|$) . \\
Similarly as in Case 1 the multiplier is bounded by $ c \frac{|\xi_2|}{N} $ 
and 
the integral by
\begin{eqnarray*}
\lefteqn{ \frac{c}{N} \int_0^{\delta} \int_* |\xi_1+\xi_2| 
\widehat{v}(\xi_1,t) 
|\xi_2| \widehat{u_2}(\xi_2,t) |\xi_3| \widehat{\overline{u_3}}(\xi_3,t) \, 
d\xi 
dt } \\
& \le & \frac{c}{N} \int_0^{\delta} |\xi_1|^{\frac{1}{2}} 
\widehat{v}(\xi_1,t) 
|\xi_2| \widehat{u_2}(\xi_2,t) |\xi_3| \widehat{\overline{u_3}}(\xi_3,t) 
|\xi_2+\xi_3|^{\frac{1}{2}} \, d\xi dt \, , 
\end{eqnarray*}  
the same bound as in Case 1, using $|\xi_1+\xi_2| \le c|\xi_1| = c 
|\xi_1|^{\frac{1}{2}} |\xi_2+\xi_3|^{\frac{1}{2}}$ . \\
\underline{Case 3:} $ |\xi_1| \sim |\xi_2| \ge N $ , $ |\xi_1+\xi_2| \le 2N 
$ . 
\\
The multiplier bound $ \frac{c}{m(\xi_1)m(\xi_2)} \le c 
\frac{|\xi_1|}{N}\frac{|\xi_2|}{N}$ implies the integral bound
\begin{eqnarray*}
\lefteqn{ cN^{-2} \int_0^{\delta} \int_* |\xi_1||\xi_2||\xi_1+\xi_2| 
\widehat{v}(\xi_1,t) \widehat{u_2}(\xi_2,t) |\xi_3| 
\widehat{\overline{u_3}}(\xi_3,t) \, d\xi dt } \\
& \le & cN^{-2} \int_0^{\delta} \int_* |\xi_1|^{\frac{1}{2}} 
|\xi_2+\xi_3|^{\frac{1}{2}} |\xi_2| 2N \widehat{v}(\xi_1,t) 
\widehat{u_2}(\xi_2,t) |\xi_3| \widehat{\overline{u_3}}(\xi_3,t) \, d\xi dt 
\\
& \le & cN^{-1} \|D_x^{\frac{1}{2}}v\|_{L^2_{xt}} \|D_x^{\frac{1}{2}}(D_x 
u_2 
D_x \overline{u_3})\|_{L^2_{xt}} \\
& \le & cN^{-1} \|v\|_{Y^{\frac{1}{2},\frac{1}{2}+}} 
\|u_2\|_{X^{1,\frac{1}{2}+}} \|u_3\|_{X^{1,\frac{1}{2}+}} \, .
\end{eqnarray*}
\underline{Case 4:} $ |\xi_1| \sim |\xi_2| \ge N $ , $ |\xi_1+\xi_2| \ge 2N$ 
. 
\\
The multiplier is bounded by
$$ \frac{m(\xi_1+\xi_2)}{m(\xi_1)m(\xi_2)} +1 \le c \frac{|\xi_1|^{1-s} 
|\xi_2|^{1-s} N^{1-s}}{N^{1-s} N^{1-s} |\xi_1+\xi_2|^{1-s}} = c 
\frac{|\xi_1|^{1-s} |\xi_2|^{1-s}}{N^{1-s} |\xi_1+\xi_2|^{1-s}} $$
and the integral by
\begin{eqnarray*}
\lefteqn{ \frac{c}{N^{1-s}}\ \int_0^{\delta} \int_* |\xi_1+\xi_2|^s 
|\xi_1|^{1-s} \widehat{v}(\xi_1,t) |\xi_2|^{1-s} \widehat{u_2}(\xi_2,t) 
|\xi_3| 
\widehat{\overline{u_3}}(\xi_3,t) \, d\xi dt } \\
& \le & \frac{c}{N} \int_0^{\delta} \int_* |\xi_1|^{\frac{1}{2}} 
\widehat{v}(\xi_1,t) |\xi_2+\xi_3|^{\frac{1}{2}} |\xi_2| 
\widehat{u_2}(\xi_2,t) 
|\xi_3| \widehat{\overline{u_3}}(\xi_3,t) \, d\xi dt \\
& \le & cN^{-1} \|D_x^{\frac{1}{2}} v\|_{L^2_{xt}} \|D_x^{\frac{1}{2}}(D_x 
u_2 
D_x \overline{u_3})\|_{L^2_{xt}} \\
& \le & cN^{-1} \|v\|_{Y^{\frac{1}{2},\frac{1}{2}+}} 
\|u_2\|_{X^{1,\frac{1}{2}+}} \|u_3\|_{X^{1,\frac{1}{2}+}}
\end{eqnarray*}
using $ |\xi_1+\xi_2|^s |\xi_1|^{1-s} \le c|\xi_1| = c|\xi_1|^{\frac{1}{2}} 
|\xi_2+\xi_3|^{\frac{1}{2}} $ and $ |\xi_2|^{1-s} \le |\xi_2| N^{-s} $ . \\
\underline{Estimate of $I_2$:} We have to show
\begin{eqnarray*}
& & \int_0^{\delta} \int_* 
\frac{m(\xi_1+\xi_2)-m(\xi_1)m(\xi_2)}{m(\xi_1)m(\xi_2)} |\xi_1+\xi_2| 
\widehat{u_1}(\xi_1,t) \widehat{\overline{u_2}}(\xi_2,t) 
\widehat{u_3}(\xi_3,t) 
\widehat{\overline{u_4}}(\xi_4,t) \, d\xi dt \\
& & \le c N^{-1} \prod_{i=1}^4 \|u_i\|_{X^{1,\frac{1}{2}+}} \, .
\end{eqnarray*}
The multiplier is bounded by $ c|\xi_1||\xi_2| N^{-2}$ , if $|\xi_1|,|\xi_2| 
\ge 
N$ , and the integral by
$$c N^{-2} \int_0^{\delta} \int_* |\xi_1| \widehat{u_1}(\xi_1,t) |\xi_2| 
\widehat{\overline{u_2}}(\xi_2,t)(|\xi_3|\widehat{u_3}(\xi_3,t) 
\widehat{\overline{u_4}}(\xi_4,t) + \widehat{u_3}(\xi_3,t) |\xi_4|
\widehat{\overline{u_4}}(\xi_4,t))  d\xi dt $$
using $|\xi_1+\xi_2|=|\xi_3+\xi_4| \le |\xi_3| + |\xi_4|$ . Strichartz' 
estimate 
gives the bound
\begin{eqnarray*}
& & c N^{-2} \|D_xu_1\|_{L^4_{xt}} \|D_xu_2\|_{L^4_{xt}} 
(\|D_xu_3\|_{L^4_{xt}} 
\|u_4\|_{L^4_{xt}} + \|u_3\|_{L^4_{xt}} \|D_xu_4\|_{L^4_{xt}}) \\
& & \le cN^{-2} \prod_{i=1}^4 \|u_i\|_{X^{1,\frac{1}{2}+}} \, .
\end{eqnarray*}
If however $|\xi_1| \ge N$ , $ |\xi_2| \le N $ , the multiplier bound 
$c|\xi_1|N^{-1}$ similarly gives the bound $cN^{-1} \prod_{i=1}^4 
\|u_i\|_{X^{1,\frac{1}{2}+}}$ . \\
\underline{Estimate of $I_3$:} It suffices to show
\begin{eqnarray*}
& & \int_0^{\delta} \int_* 
\frac{m(\xi_1+\xi_2)-m(\xi_1)m(\xi_2)}{m(\xi_1)m(\xi_2)} 
\widehat{v_1}(\xi_1,t) \widehat{u_2}(\xi_2,t) \widehat{v_3}(\xi_3,t) 
\widehat{\overline{u_4}}(\xi_4,t) \, d\xi dt \\
& & \le c N^{-1} \|v_1\|_{Y^{\frac{1}{2},\frac{1}{2}+}} 
\|u_2\|_{X^{1,\frac{1}{2}+}} \|v_3\|_{Y^{\frac{1}{2},\frac{1}{2}+}} 
\|u_4\|_{X^{1,\frac{1}{2}+}} \, .
\end{eqnarray*}
\underline{Case 1:} $ |\xi_1| \ge N$ , $ |\xi_2|\ge N $ . \\
The multiplier bound
$$ c\left(\frac{|\xi_1|}{N}\right)^{\frac{3}{4}} 
\left(\frac{|\xi_2|}{N}\right)^{\frac{3}{4}} = c \frac{|\xi_1|^{\frac{1}{2}} 
|\xi_2+\xi_3+\xi_4|^{\frac{1}{4}} |\xi_2|^{\frac{3}{4}}}{N^{\frac{3}{4}} 
N^{\frac{3}{4}}} $$
allows to estimate the integral by
$$ cN^{-\frac{3}{2}} \int_0^{\delta} \int_* |\xi_1|^{\frac{1}{2}} 
\widehat{v_1}(\xi_1,t) \langle \xi_2 \rangle \widehat{u_2}(\xi_2,t) \langle 
\xi_3 \rangle^{\frac{1}{4}} \widehat{v_3}(\xi_3,t) \langle \xi_4 
\rangle^{\frac{1}{4}} \widehat{\overline{u_4}}(\xi_4,t) \, d\xi dt \, .$$
Using H\"older's inequality with exponent 4 and Strichartz' estimate easily 
gives the desired bound.\\
\underline{Case 2:} $ |\xi_1| \ge N $ , $ |\xi_2| \le N $ (or similarly 
$|\xi_1| 
\le N$ , $ |\xi_2| \ge N $). \\
The multiplier bound $ c|\xi_1|N^{-1} \le 
c|\xi_1|^{\frac{1}{2}}|\xi_2+\xi_3+\xi_4|^{\frac{1}{2}}N^{-1}$ allows to 
estimate the integral by 
$$ cN^{-1} \int_0^{\delta} \int_* |\xi_1|^{\frac{1}{2}} 
\widehat{v_1}(\xi_1,t) 
\langle \xi_2 \rangle^{\frac{1}{2}} \widehat{u_2}(\xi_2,t) \langle \xi_3 
\rangle^{\frac{1}{2}} \widehat{v_3}(\xi_3,t) \langle \xi_4 
\rangle^{\frac{1}{2}} 
\widehat{\overline{u_4}}(\xi_4,t) \, d\xi dt \, . $$
Similarly as in Case 1 this gives the desired estimate.

Concerning the estimate for $L$ we remark that the first term on the right 
hand 
side of (\ref{CL9}) can be handled like $I_1$ and the second term like $I_4$ 
(with one derivative less). This completes the proof.
\section{Global existence}
One easily checks
$$ \|I_N u\|_{\dot{H}^1} \le c N^{1-s} \|u\|_{\dot{H}^s} $$
and also for $ 0 < s \le \frac{1}{2} $ :
$$ \|I_N v\|_{L^2} \le c N^{\frac{1}{2}-s} \|v\|_{H^{s-\frac{1}{2}}} \, .$$
Trivially one has
$$ \|I_N u\|_{L^2} \le c \|u\|_{L^2} $$
and also
$$ \|I_N u \|_{L^4} \le c \|u\|_{L^4} $$
by Mikhlin's multiplier theorem. This implies immediately for $1>s \ge 
\frac{1}{2}$ 
:
\begin{eqnarray*}
|E(I_N u,I_N v)| & \le & c(\|I_N u_x\|^2 + \|D_x^{\frac{1}{2}} I_N v\|^2 + 
\|I_N 
v\|_{L^2} \|I_N u\|_{L^4}^2) \\
& \le & c\left[N^{2(1-s)} (\|u\|_{\dot{H}^s}^2 + 
\|v\|_{\dot{H}^{s-\frac{1}{2}}}^2) + \|v\|_{L^2} \|u\|_{L^4}^2\right]
\end{eqnarray*}
and
\begin{eqnarray*}
|L(I_N u,I_N v)| & \le & c(\|I_N v\|^2 + \|I_N u\| \|I_N u_x\|) \\
& \le & c(\|v\|^2 + \|u\| N^{1-s} \|u\|_{\dot{H^s}}) \, .
\end{eqnarray*} 
Similarly, for $0< s \le \frac{1}{2}$ we get
$$ |E(I_N u,I_N v)| \le c\left[N^{2(1-s)}(\|u\|^2_{\dot{H}^s} + 
\|v\|^2_{H^{s-\frac{1}{2}}}) + N^{\frac{1}{2}-s} \|v\|_{H^{s-\frac{1}{2}}} 
\|u\|_{L^4}^2\right] $$
and
$$ |L(I_N u,I_N v)| \le c (N^{2(\frac{1}{2}-s)} \|v\|_{H^{s-\frac{1}{2}}}^2 
+ 
\|u\| N^{1-s} \|u\|_{\dot{H}^s}) \, . $$
We note that our system has a scaling invariance, i.e. if $(u,v)$ is a 
solution, 
then also
$$ u^{(\lambda)}(x,t):= \lambda^{-\frac{3}{2}} 
u(\frac{x}{\lambda},\frac{t}{\lambda^2}) \quad , \quad v^{(\lambda)}(x,t):= 
\lambda^{-2} v(\frac{x}{\lambda},\frac{t}{\lambda^2}) $$
for $\lambda > 0 $ , as one easily checks. Then
$$ \|u_0^{(\lambda)}\|_{\dot{H}^s} = \lambda^{-\frac{3}{2}} 
\|u_0(\frac{x}{\lambda})\|_{\dot{H}^s} = c \lambda^{-(s+1)} 
\|u_0\|_{\dot{H}^s} 
$$
and
$$ \|v_0^{(\lambda)}\|_{\dot{H}^{s-\frac{1}{2}}} = \lambda^{-2} 
\|v_0(\frac{x}{\lambda})\|_{\dot{H}^{s-\frac{1}{2}}} = c \lambda^{-(s+1)} 
\|v_0\|_{\dot{H}^{s-\frac{1}{2}}} \quad (\mbox{for} \, s \ge \frac{1}{2}) $$
as well as
$$ \|u_0^{(\lambda)}\|_{L^4} = c \lambda^{-\frac{5}{4}} \|u_0\|_{L^4} \; , 
\; 
\|v_0^{(\lambda)}\|_{L^2} = c \lambda^{-\frac{3}{2}} \|v_0\|_{L^2} \; , \; 
\|u_0^{(\lambda)}\|_{L^2} = c \lambda^{-1} \|u_0\|_{L^2} \, . $$
We also need
\begin{lemma}
\label{Lemma1}
For $ s \le \frac{1}{2} $ and $ \lambda \ge 1 $ the following estimate 
holds:
$$ \|v_0^{(\lambda)}\|_{H^{s-\frac{1}{2}}} \le c \lambda^{-(s+1)} 
\|v_0\|_{H^{s-\frac{1}{2}}} \, .$$
\end{lemma}
{\bf Proof:}
\begin{eqnarray*}
\|v_0^{(\lambda)}\|_{H^{s-\frac{1}{2}}} & = & \lambda^{-2} 
\|v_0(\frac{x}{\lambda})\|_{H^{s-\frac{1}{2}}} = \lambda^{-2} \| \langle \xi 
\rangle^{s-\frac{1}{2}} \widehat{v_0(\frac{x}{\lambda})}\|_{L^2} \\
& = & \lambda^{-1} \| \langle \xi \rangle^{s-\frac{1}{2}} 
\widehat{v_0}(\lambda 
\xi)\|_{L^2} = \lambda^{-\frac{3}{2}} \left(\int \left| \langle 
\frac{\eta}{\lambda}\rangle^{s-\frac{1}{2}} \widehat{v_0}(\eta)\right|^2 
d\eta 
\right)^{\frac{1}{2}} \\
& \le & c \lambda^{-\frac{3}{2}} \left[ \left(\int_{|\eta|\le 1} 
|\widehat{v_0}(\eta)|^2 d\eta\right)^{\frac{1}{2}} + \left( \int_{|\eta|\ge 
1} 
\left| \left|\frac{\eta}{\lambda}\right|^{s-\frac{1}{2}} 
\widehat{v_0}(\eta)\right|^2 d\eta\right)^{\frac{1}{2}} \right] \\
& \le & c \lambda^{-\frac{3}{2}} (1+\lambda^{-(s-\frac{1}{2})}) 
\|v_0\|_{H^{s-\frac{1}{2}}} \; \le \; c \lambda^{-(s+1)} 
\|v_0\|_{H^{s-\frac{1}{2}}} 
\, .
\end{eqnarray*}

This implies the following bounds for the modified functionals $E$ and $L$ 
for $\lambda \ge 1$ : 
\\
a) In the case $1 \ge s \ge \frac{1}{2} $ we get
\begin{eqnarray*}
|E(I_Nu_0^{(\lambda)},I_N v_0^{(\lambda)})| & \le & 
c\left(N^{2(1-s)}(\|u_0^{(\lambda)}\|_{\dot{H}^s}^2 + 
\|v_0^{(\lambda)}\|_{\dot{H}^{s-\frac{1}{2}}}^2) + \|v_0^{(\lambda)}\|_{L^2} 
\|u_0^{(\lambda)}\|_{L^4}^2 \right) \\
 & \hspace{-1.5cm} \le & \hspace{-0.9cm} c \left(N^{2(1-s)} 
\lambda^{-2(s+1)}(\|u_0\|_{\dot{H}^s} + 
\|v_0\|_{\dot{H}^{s-\frac{1}{2}}})^2 + \lambda^{-4} \|u_0\|_{L^4}^2 
\|v_0\|^2 
\right) \\
\end{eqnarray*}
Thus
$$ |E(I_N u_0^{(\lambda)},I_N v_0^{(\lambda)})| \le c_0^2 N^{2(1-s)} 
\lambda^{-2(s+1)} (1+\|u_0\|_{H^s} + \|v_0\|_{H^{s-\frac{1}{2}}})^4 \, .$$
Similarly
\begin{eqnarray*}
|L(I_Nu_0^{(\lambda)},I_Nv_0^{(\lambda)})| & \le & c(\lambda^{-3} \|v_0\|^2 
+ \lambda^{-1} \|u_0\| N^{1-s} \lambda^{-(s+1)} \|u_0\|_{\dot{H}^s}) \\
& \le & c N^{1-s} \lambda^{-(s+1)}(\|v_0\|^2 + \|u_0\| \|u_0\|_{\dot{H}^s}) 
\\
& \le & c_0 N^{1-s} \lambda^{-(s+1)} (1+\|u_0\|_{H^s} + 
\|v_0\|_{H^{s-\frac{1}{2}}})^2
\end{eqnarray*}
b) In the case $\frac{1}{4} \le s < \frac{1}{2} $ we get by Lemma 
\ref{Lemma1} :
\begin{eqnarray*}
 & & \hspace{-0.5cm} |E(I_N u_0^{(\lambda)},I_N v_0^{(\lambda)})|  \\
& & \hspace{-0.5cm}\le  c [N^{2(1-s)} \lambda^{-2(s+1)} (\|u_0\|_{\dot{H}^s} 
+ \|v_0\|_{H^{s-\frac{1}{2}}})^2 + N^{\frac{1}{2}-s} \lambda^{-(s+1)} 
\|v_0\|_{H^{s-\frac{1}{2}}} \lambda^{-\frac{5}{2}} \|u_0\|_{L^4}^2] \\
& & \hspace{-0.5cm}\le  c_0^2 N^{2(1-s)} \lambda^{-2(s+1)} (1+\|u_0\|_{H^s} 
+ \|v_0\|_{H^{s-\frac{1}{2}}})^4
\end{eqnarray*}
using the embedding $H^s \subset L^4$ for $ s \ge \frac{1}{4} $ . \\
Moreover we crudely estimate
\begin{eqnarray*}
\lefteqn{ |L(I_Nu_0^{(\lambda)},I_Nv_0^{(\lambda)})| } \\
& \le & c(N^{2(\frac{1}{2}-s)} \lambda^{-2(s+1)} 
\|v_0\|_{H^{s-\frac{1}{2}}}^2 + \lambda^{-1} \|u_0\| N^{1-s} 
\lambda^{-(s+1)} \|u_0\|_{\dot{H}^s}) \\
& \le & c_0 N^{1-s} \lambda^{-(s+1)} (1+ \|v_0\|_{H^{s-\frac{1}{2}}} + 
\|u_0\|_{H^s})^2 \, .
\end{eqnarray*} 
Now assume $N >> 1$ is given (to be chosen later), we choose $\lambda = 
\lambda(N,\|u_0\|_{H^s},$ $ \|v_0\|_{H^{s-\frac{1}{2}}})$ as follows:
$$ \lambda = N^{\frac{1-s}{1+s}} (4c_0)^{\frac{1}{1+s}} 
(1+\|v_0\|_{H^{s-\frac{1}{2}}}+\|u_0\|_{H^s})^{\frac{2}{s+1}} \, ,$$
so that $|E(I_N u_0^{(\lambda)},I_N v_0^{(\lambda)})| \le \frac{1}{4} $ and 
$|L(I_N u_0^{(\lambda)},I_N v_0^{(\lambda)})| \le \frac{1}{4} $ . Such a 
bound implies by (\ref{CL7}) the following estimate for the scaled initial 
data:
\begin{eqnarray}
\nonumber
\lefteqn{\|I_Nu_0^{(\lambda)}\|_{H^1}^2 + \|I_N 
v_0^{(\lambda)}\|_{H^{\frac{1}{2}}}^2 } \\
\label{***}
 & \le & \overline{c}^2 (|E(I_N 
u_0^{(\lambda)},I_N v_0^{(\lambda)})| + |L(I_N u_0^{(\lambda)},I_N 
v_0^{(\lambda)})|^{\frac{4}{3}} + \|I_Nu_0^{(\lambda)}\|_{L^2}^4 +1) \\
& \le & \overline{c}^2 (2+\|u_0\|_{L^2}^4) \, ,
\nonumber
\end{eqnarray}
because $ \|I_N u_0^{(\lambda)}\|_{L^2} \le \|u_0^{(\lambda)}\|_{L^2} = 
\lambda^{-1} \|u_0\|_{L^2} \le \|u_0\|_{L^2} $ for $\lambda \ge 1$ . Thus 
the local existence theorem Proposition {\ref{Proposition3} gives a solution 
on a time interval of length $\delta = \delta(\|u_0\|_{L^2})$ and
\begin{equation}
\label{****}
 \|I_N 
u^{(\lambda)} \|_{X^{1,b}_{\delta}} + \|I_N v^{(\lambda)} 
\|_{Y^{\frac{1}{2},b}_{\delta}} \le \widehat{c} \, \overline{c} 
(2+\|u_0\|_{L^2}^4)^{\frac{1}{2}} \, . 
\end{equation}
Thus Proposition \ref{Proposition4} shows
\begin{eqnarray*}
\lefteqn{|E(I_Nu^{(\lambda)}(\delta),I_N v^{(\lambda)}(\delta))| + 
|L(I_Nu^{(\lambda)}(\delta),I_Nv^{(\lambda)}(\delta))|} \\
&\le & CN^{-1} + |E(I_N u_0^{(\lambda)},I_N v_0^{(\lambda)})| + |L(I_N 
u_0^{(\lambda)},I_N v_0^{(\lambda)})| 
\end{eqnarray*}
where $C=C(\|u_0\|_{L^2})$ . We choose $N$ large enough, so that
$$ |E(I_Nu^{(\lambda)}(\delta),I_N v^{(\lambda)}(\delta))| + 
|L(I_Nu^{(\lambda)}(\delta),I_Nv^{(\lambda)}(\delta))| < 1 \, $$
and such that we can reapply the local existence theorem with time intervals 
of equal 
length (remark that $\|u_0\|_{L^2}$ is conserved) $N^{1-}$ times before the 
size of \\ 
$|E(I_Nu^{(\lambda)}(\delta),I_N v^{(\lambda)}(\delta))| + 
|L(I_Nu^{(\lambda)}(\delta),I_Nv^{(\lambda)}(\delta))|$ reaches $1$ . During 
the whole iteration process the bounds for the iterated solutions on the 
right hand side of (\ref{***}) and (\ref{****}) can be chosen uniformly.
Now, given any finite time $T$ we are able to get a solution in this way on 
$[0,T] \, ,$ provided $ N^{1-} \delta \lambda^{-2} = T $ , taking the 
scaling into account. Using the 
definition of $\lambda$ above, this means that $ N^{1-} \delta 
N^{-\frac{2(1-s)}{1+s}} = T \, . $ This can be fulfilled for a sufficiently 
large $N$ , provided $ 1 > \frac{2(1-s)}{1+s} $ $ \Longleftrightarrow 1+s > 
2(1-s) $ $ \Longleftrightarrow s > \frac{1}{3} $ . 

Thus we have proven the following global existence result:
\begin{theorem}
\label{Theorem3}
For $ 1 > s > \frac{1}{3} $ and $ (u_0,v_0) \in H^s \times H^{s-\frac{1}{2}} 
$ there exists $ b > \frac{1}{2} $ such that the Cauchy problem 
(\ref{1.1}),(\ref{1.2}),(\ref{1.3}) has a unique global solution $(u,v) \in 
X^{s,b} \times Y^{s-\frac{1}{2},b} $ with $ (u,v) \in C^0_{loc}({\bf 
R}^+,H^s \times H^{s-\frac{1}{2}}) $ .
\end{theorem}

It is also possible to show that this global solution grows at most 
polynomially in $t$. The procedure above namely shows
$$ |E( I_N u^{(\lambda)}(N^{1-}\delta),I_N v^{(\lambda)}(N^{1-} \delta))| +
|E(I_N u^{(\lambda)}(N^{1-} \delta),I_N v^{(\lambda)}(N^{1-} \delta))| \le 1 
\, .$$
This implies by (\ref{CL7}):
\begin{eqnarray*}
& & \|I_N u^{(\lambda)}(N^{1-}\delta)\|_{H^1}^2 + \|I_N 
v^{(\lambda)}(N^{1-}\delta)\|_{H^{\frac{1}{2}}}^2  \le  c(1+ \|I_N 
u^{(\lambda)}(N^{1-}\delta)\|_{L^2}^4)  \\
& & \le c(1+ \|u^{(\lambda)}(N^{1-}\delta)\|_{L^2}^4)  \le  
c(1+\|u_0^{(\lambda)}\|_{L^2}^4) \\
& & \le c(1 + \lambda^{-4} \|u_0\|_{L^2}^4)  \le  c (1+\|u_0\|_{L^2}^4)
\end{eqnarray*}
for $\lambda \ge 1$ . Thus we get
$$ \|u^{(\lambda)}(N^{1-} \delta)\|_{H^s} + \|v^{(\lambda)}(N^{1-} 
\delta)\|_{H^{s-\frac{1}{2}}} \le c \, . $$
But now
$$ \|u^{(\lambda)}(N^{1-} \delta)\|_{H^s} = \lambda^{-\frac{3}{2}} 
\|u(\frac{x}{\lambda},T)\|_{H^s} \ge c \lambda^{-(1+s)} \|u(T)\|_{H^s} $$
and similarly for $ s \ge \frac{1}{2} $ :
$$ \|v^{(\lambda)}(N^{1-} \delta)\|_{H^{s-\frac{1}{2}}} = \lambda^{-2} 
\|v(\frac{x}{\lambda},T)\|_{H^{s-\frac{1}{2}}} \ge c \lambda^{-(1+s)} 
\|v(T)\|_{H^{s-\frac{1}{2}}} \, , $$
whereas for $s < \frac{1}{2}$ :
$$ \|v^{(\lambda)}(N^{1-} \delta)\|_{H^{s-\frac{1}{2}}} \ge c 
\lambda^{-\frac{3}{2}} \|v(T)\|_{H^{s-\frac{1}{2}}} \, . $$
Because $ \lambda \sim N^{\frac{1-s}{1+s}} $ and $ T \sim N^{1-} 
\lambda^{-2} \sim N^{1-} N^{-\frac{2(1-s)}{1+s}} = N^{\frac{3s-1}{s+1}-} $ , 
thus $ N \sim T^{\frac{s+1}{3s-1}+} $ , we get $ \lambda \sim 
T^{\frac{1-s}{3s-1}+} $ , so that we have proven
\begin{theorem}
\label{Theorem4}
The global solution of Theorem \ref{Theorem3} fulfills for $t \in {\bf R}$ :
$$ \|u(t)\|_{H^s} \quad \le \quad c ( 1 + t^{\frac{(s+1)(1-s)}{3s-1}+} ) 
\quad \mbox{for} \quad 1 > s > \frac{1}{3} $$
and
\begin{eqnarray*}
\|v(t)\|_{H^{s-\frac{1}{2}}} & \le &  c ( 1 + t^{\frac{(s+1)(1-s)}{3s-1}+} ) 
\quad \mbox{for} \quad 1 > s \ge \frac{1}{2} \, , \\
 \|v(t)\|_{H^{s-\frac{1}{2}}} & \le & c ( 1 + 
t^{\frac{\frac{3}{2}(1-s)}{3s-1}+} ) \qquad \, \mbox{for} \quad \frac{1}{2} 
> s > \frac{1}{3} \, .
\end{eqnarray*}
\end{theorem}
 
\end{document}